\documentclass[11pt]{article}

\usepackage{geometry}
\date{}

\usepackage{amsthm}

\usepackage{graphics} 
\usepackage{epsfig} 
\usepackage{times} 
\usepackage{amsmath} 
\usepackage{amssymb}  
\usepackage{color}

\newtheorem{assume}{Assumption}
\newtheorem{lemma}{Lemma}
\newtheorem{prop}{Proposition}
\newtheorem{remark}{Remark}
\newtheorem{definition}{Definition}
\newtheorem{problem}{Problem}
\newtheorem{theorem}{Theorem}
\newtheorem{corollary}{Corollary}
\usepackage{algorithm}
\usepackage{algpseudocode}
\usepackage{bbm}
 \usepackage{url}

\usepackage[round]{natbib}

\title{\LARGE \bf
Learning Personalized Models with Clustered System Identification
}

\author{Leonardo F. Toso, Han Wang, and James Anderson
\thanks{This material is based upon work supported in part by DoE under grant DE-SC0022234 and NSF awards 2144634 \& 2231350. The authors are with the Department of Electrical Engineering,  Columbia University in the City of New York, New York, NY, 10027, USA. Email: \texttt{\{lt2879, hw2786, james.anderson\}@columbia.edu}.} 
}

\begin{document}

\maketitle
\thispagestyle{empty}
\pagestyle{empty}

\begin{abstract}
We address the problem of learning linear system models from observing multiple trajectories from different system dynamics. This framework encompasses a collaborative scenario where several systems seeking to estimate their dynamics are partitioned into clusters according to their system similarity. Thus, the systems within the same cluster can benefit from the observations made by the others. Considering this framework, we present an algorithm where each system alternately estimates its cluster identity and performs an estimation of its dynamics. This is then aggregated to update the model of each cluster. We show that under mild assumptions, our algorithm correctly estimates the cluster identities and achieves an approximate sample complexity that scales inversely with the number of systems in the cluster, thus facilitating a more efficient and personalized system identification process.
\end{abstract}

\section{Introduction}


System identification is the data-driven process of estimating a dynamic model of a system based on observations of the system trajectories. It plays a crucial role in aiding our understanding of complex systems and is a fundamental problem in numerous fields, including time-series analysis, control theory, robotics, and reinforcement learning \citep{aastrom1971system,ljung1998system}. The effective utilization of available data is pivotal in obtaining an accurate model estimate with a measure of uncertainty quantification. Traditional system identification, methods \citep{ljung1998system} have focused on asymptotic analysis, which, although insightful, is restrictive when dealing with small to medium sized data sets. Motivated by this, and the fact that data generation is often costly and time consuming, modern approaches focus on developing sample complexity bounds (i.e., non-asymptotic convergence analysis).

Results on the estimation of both fully~\citep{sarkar2019near,dean2020sample,simchowitz2018learning} and partially~\citep{oymak2019non,sun2020finite,tu2017non,simchowitz2019learning,zheng2020non} observed LTI systems have demonstrated that a more precise characterization of error bounds is essential for designing efficient and robust control systems \citep{dean2020sample, tu2017non, khalil1996robust}. These studies provide non-asymptotic bounds that are functions of the number of observed trajectories (see Table 1 of \citep{zheng2020non} for a summary of the bounds).

A recent body of work has begun to formalize methods for improving sample efficiency by considering data (or models generated from data) from multiple systems \citep{xin2022identifying, xin2023learning, wang2022fedsysid, CheOPD23, zhang2022multi, zhang2023meta, wang2023model}. Leveraging data from similar systems provides a promising approach although clarifying the effect of the heterogeneity in the systems and their environments is crucial. The aforementioned work have demonstrated that the benefit of collaboration typically reduces the sample complexity by a factor of the number of collaborators, when compared to the single-agent setting where each system estimate its dynamics from its own observations.


However, the approaches discussed in \citep{xin2022identifying,xin2023learning,wang2022fedsysid} compute a common estimation for all participants, thereby restricting the ability to obtain personalized estimations. Furthermore, the sample complexity bounds achieved in those studies are subject to an unavoidable heterogeneity bias that cannot be controlled by the number of trajectories or systems, thus leading to an estimation error that scales with the measure of heterogeneity among the considered systems. Specifically in~\citep{xin2022identifying,xin2023learning,wang2022fedsysid} the error of the system identification process is shown to be of order $\mathcal O(\frac{1}{\sqrt N}+\epsilon_{\text{het}})$ where $\epsilon_{\text{het}}$ characterizes the worst case heterogeneity and $N$ is the number of trajectories across all systems.

Personalization in collaborative settings aims to provide tailored solutions (e.g. model estimates) to individual agents with distinct objectives, while enabling inter-agent collaboration (e.g. model sharing). This encompasses diverse topics such as representation learning \citep{zhang2023meta, bengio2013representation,finn2017model, collins2021exploiting} and clustering \citep{xu2005survey}, both widely studied in machine learning and data analysis. The present work address the aforementioned challenges by leveraging clustering techniques to achieve personalized model estimations. The idea is simple: cluster systems into groups that have identical system dynamics, and then apply collaborative learning algorithms to the clusters in order to improve sample complexity (by reducing the heterogeneity induced error $\epsilon_{\text{het}}$) and achieve personalization even for heterogeneous settings.

Recent work on clustered federated learning that includes \citep{ghosh2020efficient}, \citep{ghosh2022improved}, \citep{sattler2020clustered} have shown the potential of clustering techniques to collaboratively train models in heterogeneous settings with non-i.i.d. data. Building upon this success, this paper aims to apply clustering to the system identification problem, which poses unique challenges due to the dynamical nature of the system that results in non-i.i.d. data. This is in contrast to the linear regression and model training settings explored in the aforementioned work. Further details on these challenges are discussed later.

Specifically, we investigate the scenario where we have $M$ dynamical systems, with each of them belonging to one of $K$ different system types (which we refer to as a ``cluster''). Which cluster a system belongs to is not initially disclosed. Our objective is to simultaneously identify the correct cluster identities for each of the $M$ systems and obtain a system model by collaboratively learning with the systems in the same cluster. Our approach can lead to significant reductions in the amount of data required to accurately estimate the system models, as illustrated in the following theorem.

\begin{theorem} (main result, informal) Suppose the $K$ system types are sufficiently different, and we observe the same number of trajectories from each system. Then, for a given cluster, with high probability, the estimation error between the learned and ground truth model is bounded by:
\begin{align*}
  \begin{array}{c}  \text{estimation} \text{ error} \end{array}  &\lesssim \frac{\text{1}}{\sqrt{\text{\# systems} \times \text{\# trajectories}}}+ \begin{array}{c}  \text{misclass.}  \text{ rate}\end{array} ,
\end{align*}
with
\[
\begin{array}{c} \text{misclass.} \text{ rate}  \end{array} \lesssim \exp(-\text{\#~trajectories}\times \text{misclass. const.}).
\]
\noindent where \text{\#systems} denotes the number of systems in the cluster, and \text{\#trajectories} represents the number of trajectories observed by each of them.
\end{theorem}

The first term captures the error in learning the system dynamics from systems' observations within the same cluster. It shows what one would hope; as the number of systems and observations increase, the error decreases. However, this speedup does not come for free. The second term is the penalty paid for assigning one of the $M$ systems to one of the incorrect $K$ clusters. One of the main results from our work is to show that \emph{both terms} can be controlled by adjusting the number of observed trajectories. Moreover, the misclassification rate is dominated by the first term, thus leading to a an approximate sample complexity that is scale inversely with the number of system within the cluster. This is in stark contrast to \citep{xin2022identifying,xin2023learning,wang2022fedsysid} which is where the heterogeneity introduces a bias $\epsilon$ which is not a function of the number of systems or the volume of data at our disposal. Our work shows that by controlling both sources of error, our approach can accurately estimate the system dynamics with fewer samples, when compared to the single agent case, and provides better estimation in heterogeneous settings when compared to \citep{xin2022identifying,xin2023learning,wang2022fedsysid}.

\emph{Contributions}: This is the first work to introduce clustering in order to provide sample complexity gains to the collaborative system identification problem. We derive an upper bound on the estimation error (Theorem \ref{theorem:convergence}) that decomposes into two terms (as shown above), where each term can be controlled by adjusting the number of observed trajectories. We offer theoretical guarantees on the probability of cluster identity misclassification (Lemma \ref{lemma:prob_misclassification}) and thus convergence (Corollary \ref{corollary:convergence}). We show that under a mild assumption on the number of observed trajectories, our approach correctly estimates the cluster identities, with high probability. Moreover, we show that our method achieves an improved convergence rate when compared to the single-agent system identification process. In contrast to the federated setting~\citep{wang2022fedsysid,CheOPD23} and that of~\citep{xin2022identifying, xin2023learning}, we are able to provide personalized models as opposed to a single generic model, thus expanding the use cases for collaborative system identification.

\subsection{Notation}

Given a matrix $G \in \mathbb{R}^{m\times n}$, the Frobenius norm of $G$ is denoted by $\|G\|_F= \sqrt{Tr(GG^\top)}$.  $\|G\|=\sigma_{\text{max}}(G)$, where $\sigma_{\text{max}}(G)$ is the largest singular value of $G$. Consider a symmetric matrix $\Sigma$, $\lambda_{\text{min}}(\Sigma)$ and $\lambda_{\text{max}}(\Sigma)$ denote its minimum and maximum eigenvalues, respectively. For systems, we use superscript $(i)$ to denote the system index and subscript $t$ for time. For models, subscript denotes the cluster identity, and superscript $(r)$ is the iteration counter.

\section{Problem Formulation and Algorithm}\label{sec:problem_formulation}

Consider $M$  linear time-invariant (LTI) systems
\begin{align}\label{eq:FL_dynamical_systems}
    x^{(i)}_{t+1} = A^{(i)}x^{(i)}_t + B^{(i)}u^{(i)}_t + w^{(i)}_t, \;\ t=0,1,\ldots,T-1
\end{align}
 where $x^{(i)}_t \in \mathbb{R}^{n_x}$, $u^{(i)}_t \in \mathbb{R}^{n_u}$ and $w^{(i)}_t \in \mathbb{R}^{n_x}$ are the state, input, and process noise at time $t$, for system $i \in [M]$. We assume that $\{u^{(i)}_t\}_{t=1}^{T-1}, \{w^{(i)}_t\}_{t=1}^{T-1}$ are random vectors distributed according to $u^{(i)}_{t} \stackrel{\text{i.i.d.}}{\sim} \mathcal{N}\left(0, \sigma_{u,i}^{2} I_{n_u}\right)$ and $w^{(i)}_{t}\stackrel{\text{i.i.d.}}{\sim} \mathcal{N}\left(0, \sigma_{w,i}^{2} I_{n_x}\right)$. Furthermore, it is assumed that $x^{(i)}_{0} \stackrel{\text{i.i.d.}}{\sim} \mathcal{N}\left(0, \sigma_{x,i}^{2} I_{n_x}\right)$.  
 
 We consider the setting where we have access to $M$ datasets corresponding to observed system trajectories. Each of the datasets is generated by one of $K$ different systems. We consider the case where $K \ll M$.  We will from now on refer to the $K$  types of different systems as ``clusters'', which we label as $\mathcal{C}_1,\ldots, \mathcal{C}_K$. We denote $(A_j, B_j)$ as the ground truth system matrices of  cluster $j \in [K]$. That is, $A^{(i)} = A_j$, and $B^{(i)} = B_j$, for any $i \in \mathcal{C}_j$. Note that due to the noise in model~\eqref{eq:FL_dynamical_systems}, two datasets generated by cluster $\mathcal C_j$ will be different. 
 
The state-input pair of a single trajectory $\{x^{(i)}_t, u^{(i)}_t\}$ of system $i \in \mathcal{C}_j$ is referred to as \textit{rollout}. We consider the setting where multiple rollouts of length $T$ are collected and stored as  $\left\{x^{(i)}_{l,t}, u^{(i)}_{l,t}\right\}_{t=0}^{T-1}$, for $l=  1,\hdots  {{N_i}}$, with $l$ denoting the $l$-th rollout and $t$ the $t$-th time-step of the corresponding rollout. Thus, for any system $i \in \mathcal{C}_j$ and cluster $j \in [K]$, the system dynamics is described by:  
\begin{align}\label{eq:FL_augmented_system}
 x^{(i)}_{l,t+1} = \Theta_j z^{(i)}_{l,t} + w^{(i)}_{l,t} \quad  \forall \ 1 \le l \le {{N_i}} \text{ and } 0\le t \le T-1,   
\end{align}
where $ z^{(i),\top}_{l,t} \triangleq\left[
x^{(i),\top}_{l,t} \;\ u^{(i),\top}_{l,t} \right] \in \mathbb{R}^{n_x+n_u}$ corresponds to the augmented state-input pair of system $i \in \mathcal{C}_j$ over rollout $l$ at time $t$, and $\Theta_j \triangleq [A_j \;\ B_j]$ denotes the concatenation of the ground truth system matrices $A_j$ and $B_j$. The state update $x^{(i)}_{l,t+1}$ can be expanded recursively as follows:
$$
x^{(i)}_{l,t}=G^{(i)}_{t}\left[\begin{array}{c}
u^{(i)}_{l,0} \\
\vdots \\
u^{(i)}_{l,t-1}
\end{array}\right]+F^{(i)}_{t}\left[\begin{array}{c}
{w}^{(i)}_{l,0} \\
\vdots \\
{w}^{(i)}_{l,t-1}
\end{array}\right]+{A_j^t}{x}^{(i)}_{l,0},
$$
where, $
G^{(i)}_{t} \triangleq\left[\begin{array}{llll}
A_j^{t-1} B_j & A_j^{t-2} B_j & \cdots & B_j
\end{array}\right]$ and
$F_{t} \triangleq \left[\begin{array}{llll}
A_j^{t-1} & A_j^{t-2} & \cdots & I_{n_x}
\end{array}\right] 
$
for all $t \geq 1$.  

The state-input pair ${z}^{(i)}_{l,t}$ is distributed according to a Gaussian distribution with zero mean and covariance matrix $\Sigma^{(i)}_t $, where, 
$$
{\Sigma}^{(i)}_{0} \triangleq\left[\begin{array}{cc}
\sigma_{{x,i}}^{2} I_{n_x} & 0 \\
0 & \sigma_{{u,i}}^{2} I_{n_u}
\end{array}\right] \succ 0,\quad \text{for }t=0,
$$
and 
\begin{small}
$$
{\Sigma}^{(i)}_{t} \triangleq\left[\begin{array}{cc}
\sigma_{{u,i}}^{2} {G}^{(i)}_{t} {G}^{(i),\top}_{t}+\sigma_{{w,i}}^{2} {F}^{(i)}_{t} {F}^{(i),\top}_{t}+\sigma_{{x,i}}^{2} {A}_j^{t} ({A}_j^{t})^{\top} & 0 \\
0 & \sigma_{{u,i}}^{2} I_{n_u}
\end{array}\right],
$$
\end{small}
for all $t\geq 1$ and $i \in \mathcal{C}_j$, $\forall j \in [K]$, as detailed in \citep{wang2022fedsysid}.

Next, we define the offline batch matrices for each system $i \in \mathcal{C}_j$, $\forall j \in [K]$. For a single rollout $l$, the data is concatenated according to
${X}^{(i)}_{l}=\left[\begin{array}{lll}{x}^{(i)}_{l,T} & \cdots & {x}^{(i)}_{l,1}\end{array}\right] \in \mathbb{R}^{n_x \times T}, \quad {Z}^{(i)}_l=\left[\begin{array}{lll}{z}^{(i)}_{l,T-1}&\cdots& {z}^{(i)}_{l,0}\end{array}\right] \in \mathbb{R}^{(n_x+n_u) \times T},$ and ${W}^{(i)}_l=\left[\begin{array}{llll}{w}^{(i)}_{l,T-1} & \cdots & {w}^{(i)}_{l,0}\end{array}\right] \in \mathbb{R}^{n_x \times T}$.  This is then further stacked to construct the batch matrices ${X}^{(i)}=\left[\begin{array}{lll}{X}^{(i)}_1 & \ldots & {X}^{(i)}_{{{N_i}}}\end{array}\right] \in \mathbb{R}^{n_x \times {{N_i}} T}, \quad {Z}^{(i)}=\left[\begin{array}{lll}{Z}^{(i)}_{1} & \cdots & {Z}^{(i)}_{{{N_i}}}\end{array}\right] \in \mathbb{R}^{(n_x+n_u) \times {{N_i}} T},$ and ${W}^{(i)}=\left[\begin{array}{lll}{W}^{(i)}_{1} & \cdots & {W}^{(i)}_{{{N_i}}}\end{array}\right]\in\mathbb{R}^{n_x \times {{N_i}} T}$. Therefore, for each system $i \in \mathcal{C}_j$, $\forall j \in [K]$, its state, input, noise, and model parameters are related according to
\begin{align}
    X^{(i)} = \Theta_j Z^{(i)} + W^{(i)},
\end{align}
\vspace{-0.1em}
where each column of $Z^{(i)}$ and $W^{(i)}$ are sampled according to Gaussian distributions with zero means and covariance matrices $\Sigma^{(i)}_t$, $\sigma_{w,i}^{2} I_{n_x}$, respectively. With that said, we are now able to introduce the clustered system identification problem. 

\begin{problem}\label{problem:clustered_sysid}
We consider $M$ dynamical systems as in \eqref{eq:FL_dynamical_systems} that are equipped with batch matrices $X^{(i)}, Z^{(i)}$, and $W^{(i)}$. Each system $i \in [M]$ is associated with its own cost function $C^{(i)}(\Theta)=\|{X}^{(i)}-{{\Theta}} {Z}^{(i)}\|_{F}^{2}$, and is unaware of its cluster identity. We aim to estimate the systems' cluster identities $\widehat{\mathcal{C}}_1,\ldots, \widehat{\mathcal{C}}_K$ and use it to estimate a model $\widehat{\Theta}_j = [\widehat{A}_j\;\ \widehat{B}_j]$ which is close to the ground truth $\Theta_j$, $\forall j \in [K]$. \end{problem}

To obtain a faster and more accurate estimation, we frame the system identification problem in the setting where systems within the same cluster can leverage data from each other. Further in this paper, we provide theoretical guarantees to support these statements.

The problem described above can be framed into an alternating optimization problem, as the actual cluster identity of each system (i.e., $\mathcal{C}_1, \ldots, \mathcal{C}_K$) is not disclosed to the systems in advance. Therefore, our objective is twofold: firstly, we aim to classify the correct cluster identities of the systems by employing the Mean Square Error (MSE) as the clustering criterion, with the resulting output being the cluster estimation (CE); secondly, we use that estimation to identify the model dynamics of each cluster with a model estimation (ME) step.
 Next, we introduce our clustered system identification algorithm to solve this problem.
\begin{algorithm}
\caption{\texttt{Clustered System Identification}} \label{algorithm:clustered_fedsysid}
\begin{algorithmic}[1]
\State \textbf{Initialization:} number of clusters $K$,  step-size $\eta_j$, and model initialization $\widehat{\Theta}^{(0)}_j $ $\forall j \in [K]$,
\State \textbf{for} each iteration $r=0, 1, \ldots, R-1$ \textbf{do}
\State \quad \quad The systems receive the models $\{\widehat{\Theta}^{(r)}_1, \ldots,\widehat{\Theta}^{(r)}_K\}$, $ \forall j \in [K]$,
\State \quad \quad  \textbf{Cluster estimation (CE):}
\State \quad \quad \quad \textbf{for} each system $i \in [M]$
\State \quad \quad \quad \quad $\hat{j}=\text{argmin}_{j \in [K]} \|{X}^{(i)} -\widehat{\Theta}^{(r)}_{j} {Z}^{(i)}\|_{F}^{2}$, 
\State \quad \quad \quad \quad  define  $e_i=\left\{e_{i, j}\right\}_{j=1}^K$ with $e_{i, j}=\mathbbm{1} \{j=\widehat{j}\}$,
\State \quad \quad \quad  \textbf{end for}
\State \quad \quad  \textbf{Model estimation (ME):}
\State \quad \quad   \begin{small} $\widehat{\Theta}^{(r+1)}_j =\widehat{\Theta}^{(r)}_{j} + \frac{2\eta_j}{\sum_{i \in [M]}e_{i,j}} \sum_{i \in [M]} e_{i,j} (X^{(i)} - \widehat{\Theta}^{(r)}_{j}Z^{(i)})Z^{(i),\top}$\end{small} for all $j \in [K]$
\State \textbf{end for}
\State \textbf{Return} $\widehat{\Theta}^{(R)}_{j}$ for all $j \in [K]$.
\end{algorithmic}
\end{algorithm}

The initial step of Algorithm \ref{algorithm:clustered_fedsysid} involves the initialization of the number of clusters and the provision of an initial guess for the dynamics of each cluster. Subsequently, the algorithm iterates from line 2 to 11, during which each system estimates its corresponding cluster identity and stores this information in the form of a one-hot encoding vector denoted by $e_i$. The one-hot encoding vector comprises $K$ elements, with one in the position of the estimated cluster identity and zero elsewhere. After the estimation of the cluster identity, the cluster model is updated by performing a single gradient descent iteration in line 10, with the gradient being the average of the gradients of each individual system's cost function that belongs to the cluster.

\begin{remark}
Note that Algorithm \ref{algorithm:clustered_fedsysid} is an alternating minimization algorithm, where it performs an iterative clustering step followed by a model estimation process. Prior to the start of collaboration, each system $i \in [M]$ collects data and stores it in batch matrices $X^{(i)}, Z^{(i)}$, and $W^{(i)}$. Moreover, it is worth noting that Algorithm \ref{algorithm:clustered_fedsysid} uses the same batch matrices for both cluster identity and  model estimation.
\end{remark}

The following definitions and assumptions are required in order to analyze Algorithm 1. Subsequently, we provide the intuition behind them.

\begin{definition} \label{def:minimum_separation}
The minimum and maximum separation between the clusters are defined as 
\begin{align*}
\Delta_{\text{min}} \triangleq \min_{j \neq j^{\prime}} \| \Theta_j - \Theta_{j^\prime}\| \quad \text{and} \quad 
\Delta_{\text{max}} \triangleq \max_{j \neq j^{\prime}} \| \Theta_j - \Theta_{j^\prime}\|,
\end{align*}
 respectively.
\end{definition}

We define $\rho^{(i)} \triangleq \frac{\Delta^2_{\min}}{\sigma^2_{w,i}}$ as the signal-to-noise ratio $\forall i \in [M]$.

\begin{assume} \label{assumption:warm_initialization}
 The initial model estimate $\widehat{\Theta}_j^{(0)}$ satisfy $\|\widehat{\Theta}_j^{(0)}-{\Theta}_j\| \leq\left(\frac{1}{2}-\alpha^{(0)}\right) \Delta_{\min}, \forall j \in[K]$, where $0<\alpha^{(0)}<\frac{1}{2}$.
\end{assume}

\begin{assume} \label{assumption:Ni_nx}
 For any fixed and small $\delta$, the number of trajectories satisfies $N_i n_x \gtrsim \left(\frac{\rho^{(i)}\|\Sigma^{(i)}_t\| + \sqrt{n_x}}{\alpha^{(0)} \rho^{(i)} \|\Sigma^{(i)}_t\| }\right)^2
\log(\frac{MT}{\delta}) $, for all $i \in [M]$. We also assume that $\Delta_{\min} \gtrsim 1 + \Delta_{\max} \sum_{i \in [M]}\sum_{t=0}^{T-1} \exp \left(-c N_in_x \left(\frac{\alpha^{(0)} \rho^{(i)} \|\Sigma^{(i)}_t\|}{\rho^{(i)}\|\Sigma^{(i)}_t\| + \sqrt{n_x} }\right)^2\right)$ for some constant $c$.
\end{assume}

Assumption \ref{assumption:warm_initialization} implies that the initial guess for the model estimates is superior to a random initialization. This assumption  is standard for alternating minimization algorithms, particularly for learning mixture models \citep{balakrishnan2017statistical}. The condition on the number of trajectories in Assumption \ref{assumption:Ni_nx} is a common requirement in the concentration bound analysis. This is used to guarantee that the cluster estimation procedure of Algorithm \ref{algorithm:clustered_fedsysid} correctly estimate the cluster identities, with high probability. Note that this is a mild assumption since for well-behaved systems where $\Sigma^{(i)}_t$ is well conditioned, $N_i n_x$ is typically in the same or superior to the order of $\log\left(\frac{MT}{\delta}\right)$. The condition on $\Delta_{\min}$ in Assumption \ref{assumption:Ni_nx} is to ensure that any two clusters are well-separated. This is a standard assumption in the literature of clustering \citep{dunn1974well, kumar2010clustering}. Similar assumptions are exploited in \citep{ghosh2020efficient} in the context of the linear regression problem.

\section{Theoretical Guarantees}\label{sec:theoretical_guarantees}

We begin our analysis by examining a single iteration of Algorithm \ref{algorithm:clustered_fedsysid}. For simplicity, we omit the superscript $r$ that denotes the iteration counter. Let us assume that we have the current estimated model $\widehat{\Theta}_j$ for all clusters $j \in [K]$ at a given iteration, such that $\|\widehat{\Theta}_j-{\Theta}_j\| \leq\left(\frac{1}{2}-\alpha\right) \Delta_{\min}$ for all $j \in[K]$, with $0<\alpha<\frac{1}{2}$.
\vspace{-1em}
\subsection{Probability of Cluster Identity Misclassification}
\vspace{-0.5em}
Consider a system $i \in [M]$ within cluster $\mathcal{C}_j$.  Let $\mathcal{M}_i^{j, j^{\prime}}$ be the event in which system $i$ is inaccurately classified as belonging to cluster $\mathcal{C}_{j^\prime}$. 
The event when system $i$ is \emph{correctly} classified is denoted as $\mathcal{M}_i^{j, j}$. The following lemma provides an upper bound on the probability of misclassification.

\begin{lemma} \label{lemma:prob_misclassification} Suppose that  $i \in \mathcal{C}_j$. There exist universal constants $c_1$ and $c_2$, such that for any $j^{\prime} \neq j$,
\vspace{-0.5em}
\begin{align*}
\mathbb{P}\left\{\mathcal{M}_i^{j, j^{\prime}}\right\}
\leq c_1 \sum_{t=0}^{T-1} \exp\left(-c_2 N_in_x \left(\frac{\alpha \rho^{(i)} \|\Sigma^{(i)}_t\|}{\rho^{(i)}\|\Sigma^{(i)}_t\| + \sqrt{n_x} }\right)^2\right).
\end{align*}  
\end{lemma}

We prove Lemma \ref{lemma:prob_misclassification} in Appendix \ref{proof:probability_of_misclassification}. By combining Lemma \ref{lemma:prob_misclassification} with the condition on $N_in_x$ from Assumption~\ref{assumption:Ni_nx}, our algorithm can ensure that the probability of misclassifying system $i$ to cluster $\mathcal{C}_{j^\prime}$ is at most $\delta$, where $\delta$ can be arbitrarily small. Moreover, it is noteworthy that if we assume the data $X^{(i)}$, $Z^{(i)}$, and $W^{(i)}$ to be i.i.d. with $T=1$ and $n_x=1$, and the columns of $Z^{(i)}$ to have an identity covariance matrix, we can recover the probability of misclassification in the linear regression problem, as discussed in \citep{ghosh2020efficient}. 
\vspace{-1em}
\subsection{Convergence Analysis}
\vspace{-0.5em}
We now examine the convergence of Algorithm \ref{algorithm:clustered_fedsysid}. The theorem below is a single-iteration convergence analysis of our algorithm.
Here we assume that, at a given iteration, an estimation $\widehat{\Theta}_j$ is obtained, which closely approximates the true model $\Theta_j$, i.e., $\|\widehat{\Theta}_j-{\Theta}_j\| \leq\left(\frac{1}{2}-\alpha\right) \Delta_{\min}$, $\forall j \in[K]$ and $0<\alpha<\frac{1}{2}$. We demonstrate that $\widehat{\Theta}_j$ converges to $\Theta_j$ up to a small bias.

\begin{theorem} \label{theorem:convergence} For any fixed $0<\delta<1$, with $N_{i} \geq$ $\max \left\{8(n_x+n_u)+16 \log \frac{2MT}{\delta}, (4 n_x+2 n_u) \log \frac{MT}{\delta}\right\}$, $\forall i \in [M]$, and selected step-size $\eta_j = \frac{|{\mathcal{C}}_j|}{\lambda_{\min}\left(\sum_{i \in {\mathcal{C}}_j } N_i \sum_{t=0}^{T-1} {\Sigma}_{t}^{(i)}\right)}$, with probability at least $1-3\delta$, it holds that,
\begin{align}
    &\| \widehat{\Theta}^+_j -\Theta_j \| \leq \frac{1}{2}\|  \widehat{\Theta}_j -\Theta_j\| + \Bar{c}_0 \times \frac{1}{\sqrt{\sum_{i \in \widehat{\mathcal{C}}_j} N_i}} +\bar{c}_1\Delta_{\max} \sum_{i \in [M] }\sum_{t=0}^{T-1} \exp\left(-\bar{c}_2 N_in_x \left(\frac{\alpha \rho^{(i)} \|\Sigma^{(i)}_t\|}{\rho^{(i)}\|\Sigma^{(i)}_t\| + \sqrt{n_x} }\right)^2\right),
\end{align}
for all $j \in [K]$, where $\bar{c}_0$, $\bar{c}_1$, $\bar{c}_2 > 0$ are problem dependent constants.
\end{theorem}

The proof of Theorem \ref{theorem:convergence} is detailed in Appendix \ref{proof:convergence}. This theorem provides an upper bound for the estimation error per iteration of our algorithm. Specifically, this bound consists of three terms. The first term is a contraction term that decreases to zero as the number of iterations increases. The second term is a constant error that decreases as the total number of observed trajectories by the systems within the cluster increases. The final term is the misclassification rate, which decays exponentially with the number of observed trajectories.

Note that although our setting is different from \citep{ghosh2020efficient}, which leads to a different estimation error expression, our per-iteration estimation error is also composed of a contractive term added to a constant error that can be controlled by the amount of data (i.e., the number of observed trajectories). We proceed to show the convergence of our algorithm by demonstrating that $\alpha^{(r)}$ is non-decreasing throughout iterations and using Assumptions \ref{assumption:warm_initialization} and \ref{assumption:Ni_nx} to show that $\|\widehat{\Theta}^{(r+1)}_j - \Theta_j\| \leq \|\widehat{\Theta}^{(r)}_j - \Theta_j\|$ for all $r \in [R]$.

Therefore, equipped with the aforementioned result, the following corollary characterizes the convergence of Algorithm \ref{algorithm:clustered_fedsysid} by providing the number of iterations required to attain a certain small and near optimal error $\epsilon$, i.e.,  $\|\widehat{\Theta}^{(R)}_j - \Theta_j\| \leq \epsilon$, for all clusters $j \in [K]$.

\begin{corollary} \label{corollary:convergence}
Frame the hypotheses of Theorem~\ref{theorem:convergence} and Assumptions \ref{assumption:warm_initialization} and \ref{assumption:Ni_nx}. Select the step-size as $\eta_j = \frac{|\mathcal{C}_j|}{\lambda_{\min}\left(\sum_{i \in \mathcal{C}_j } N_i \sum_{t=0}^{T-1} {\Sigma}_{t}^{(i)}\right)}$ for all $j \in [K]$. Then, after $R\geq  2+\log(\frac{\Delta_{\min}}{4\epsilon})$ parallel iterations, we have $\|\widehat{\Theta}^{(R)}_j - \Theta_j\| \leq \epsilon$, with 
\vspace{-1em}
\begin{align}\label{eq:epsilon_corollary}
    \epsilon&= \tilde{c}_0 \times \frac{1}{\sqrt{\sum_{i \in {\mathcal{C}}_j} N_i}}  +\tilde{c}_1\Delta_{\max} \sum_{i \in [M] }\sum_{t=0}^{T-1} \exp\left(-\tilde{c}_2 N_in_x \left(\frac{\rho^{(i)} \|\Sigma^{(i)}_t\|}{\rho^{(i)}\|\Sigma^{(i)}_t\| + \sqrt{n_x} }\right)^2\right),
\end{align}
for all $j \in [K]$, where $\tilde{c}_0$, $\tilde{c}_1$, $\tilde{c}_2 > 0$ are problem dependent constants.
 \end{corollary}

The proof of Corollary 1 can be found in Appendix \ref{proof:corollary_convergence}. Our proof builds upon similar arguments as in \citep{ghosh2020efficient}, which considers the linear regression setting. To establish the non-decreasing property of $\alpha^{(r)}$ for all $r \in [R]$ and a decrease in the additive error term over the iterations, we rely on Assumptions \ref{assumption:warm_initialization} and \ref{assumption:Ni_nx}. Furthermore, we demonstrate that our algorithm achieves a sufficiently large value of $\alpha^{(r)} \geq \frac{1}{4}$ after only a small number of iterations $R\geq 2$. This indicates that after a suitable number of iterations, our Algorithm \ref{algorithm:clustered_fedsysid} produces an estimation error that scales down with the number of systems within the cluster, and is independent of the initial closeness parameter $\alpha^{(0)}$. 
 
 This corollary highlights the benefits of collaboration. It demonstrates that the estimation error scales inversely with the number of agents within a cluster, implying that as the number of systems in the cluster increases, this error decreases. This leads to a smaller error when compared to the single agent setting, where each system estimates its dynamics using \emph{only} its own observations.

Importantly, the presented error bound differs from that of \citep{wang2022fedsysid}. Here the misclassification rate exponentially decays with the number of observed trajectories, whereas the heterogeneity bias $\epsilon_{\text{het}}$ in \citep{wang2022fedsysid} cannot be controlled by the number of trajectories. This indicates that under heterogeneous settings where the systems are significantly different, our clustering-based approach outperforms \citep{wang2022fedsysid} by providing better control over the sources of error. However, it is worth mentioning that when the systems are similar and personalization is not required, the approaches introduced in~\citep{xin2022identifying,xin2023learning,wang2022fedsysid} may be more favorable as their error bounds scale down with the total number of systems and do not necessitate a clustering step.

\section{Numerical Results} \label{sec:numerical_results}

The following simulations\footnote{Code can be downloaded from \url{https://github.com/jd-anderson/cluster-sysID}}  illustrate the efficiency of Algorithm \ref{algorithm:clustered_fedsysid}. Our analysis considers $M=50$ systems, each described by an LTI model as in \eqref{eq:FL_dynamical_systems} where $K=3$ clusters and  the number of systems in each cluster is $|\mathcal{C}_1|=10$, $|\mathcal{C}_2|=24$, and $|\mathcal{C}_3|=16$. The systems matrices for each cluster are described as follows:
\begin{align*}
   & A_1=\begin{bmatrix}
        0.5 & 0.3 &  0.1\\
        0.0  &  0.2 &  0.0\\
        0.1 &   0.0 & 0.3\\
    \end{bmatrix},\;\ A_2=\begin{bmatrix}
        -0.3 &  0.0  &  0.0\\
        0.1 & 0.4 &  0.0\\
        0.2 & 0.3 & 0.5\\
    \end{bmatrix},\;\ A_3=\begin{bmatrix}
       -0.1 &  0.1 &    0.1\\
        0.1 &  0.15 &   0.1\\
        0.1 &   0.0  &   0.2
    \end{bmatrix},\\
     &B_1 = \begin{bmatrix}
        1  &   0.5\\
        0.1 &    1\\
        0.75 & 1.5
    \end{bmatrix},\;\ B_2 = \begin{bmatrix}
        1  &   0.5\\
        0.1 &    1\\
        0.75 & 1.5
    \end{bmatrix},\;\ B_3 = \begin{bmatrix}
        0.8 &  0.1\\
        0.1 &  1.5\\
        0.4 & 0.8
    \end{bmatrix},
\end{align*}
where the initial state, input, and process noise standard deviations, for each cluster, are set to $\sigma_{x,i}=\sigma_{u,i}=\sigma_{w,i}=0.11$, $ \forall i \in \mathcal{C}_1$, $\sigma_{x,i}=\sigma_{u,i}=\sigma_{w,i}=0.12$,  $\forall i \in \mathcal{C}_2$, and $\sigma_{x,i}=\sigma_{u,i}=\sigma_{w,i}=0.05$, $\forall i \in \mathcal{C}_3$. We consider the same number of trajectories $N_i = 100$ for all $i \in [M]$. Moreover, the trajectory length is set to $T=50$. We use a fixed step-size $\eta_j=10^{-3}$, $\forall j \in [3]$. For each iteration $r$, the estimation error $e^{(j)}_r$ is defined as the spectral norm distance between the estimated model $\widehat{\Theta}^{(r)}_j$ and the ground truth model $\Theta_j$, i.e., $e^{(j)}_r = \| \widehat{\Theta}^{(r)}_j - \Theta_j\|$, for all clusters $j \in [K]$.

\begin{figure*}
    \centering
    \includegraphics[width=0.3\textwidth]{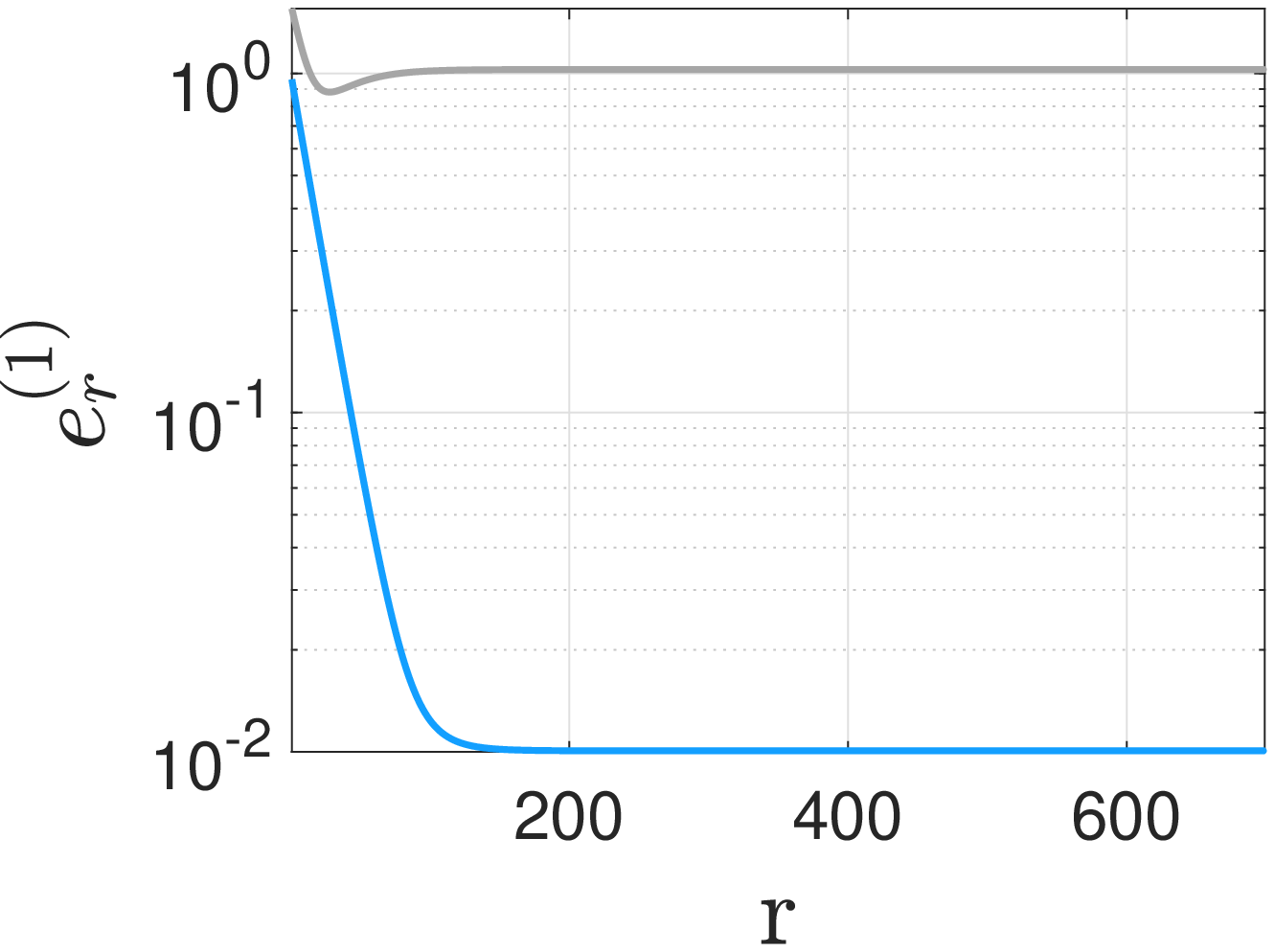}\label{fig:with_without_cluster_1}
    \includegraphics[width=0.3\textwidth]{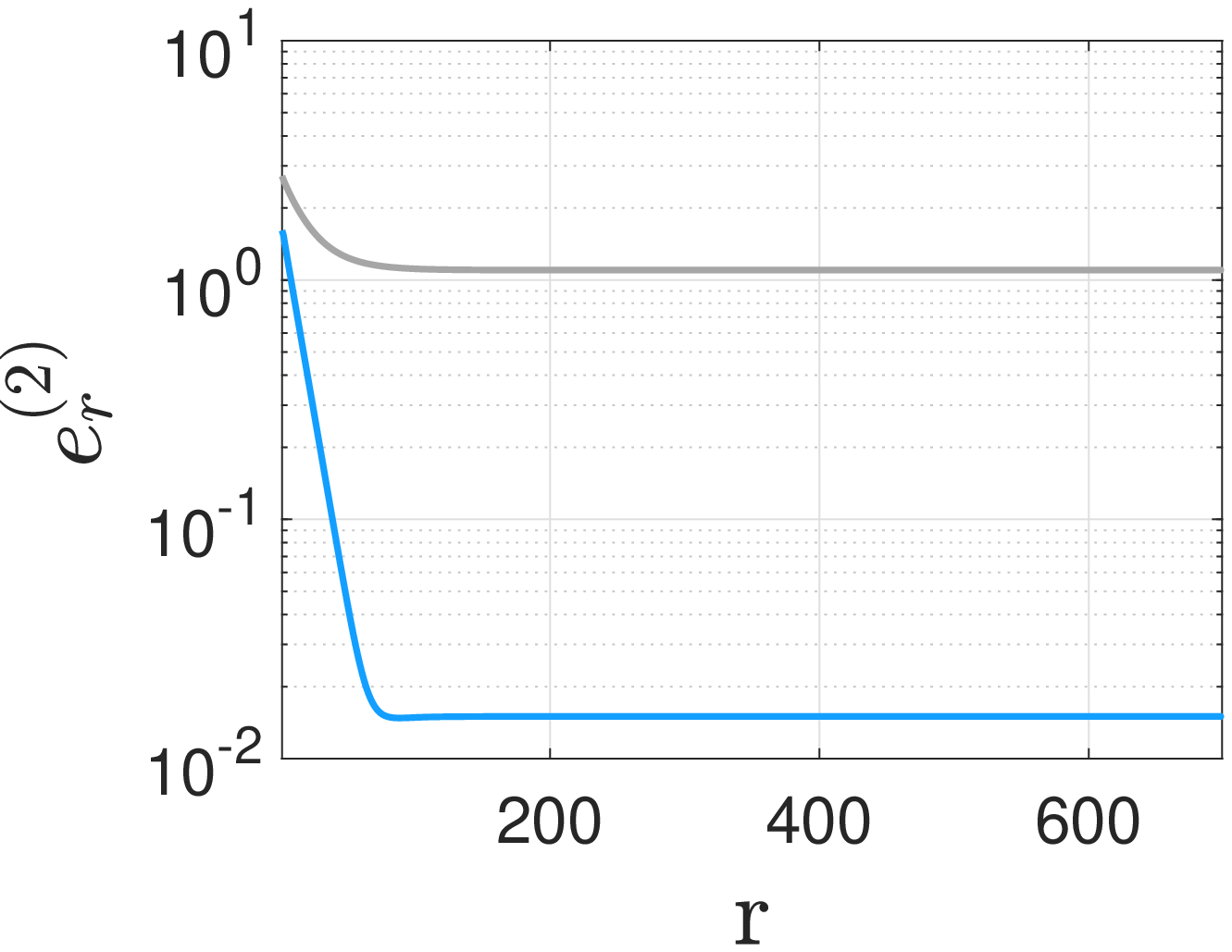}\label{fig:with_without_cluster_2}
   \includegraphics[width=0.3\textwidth]{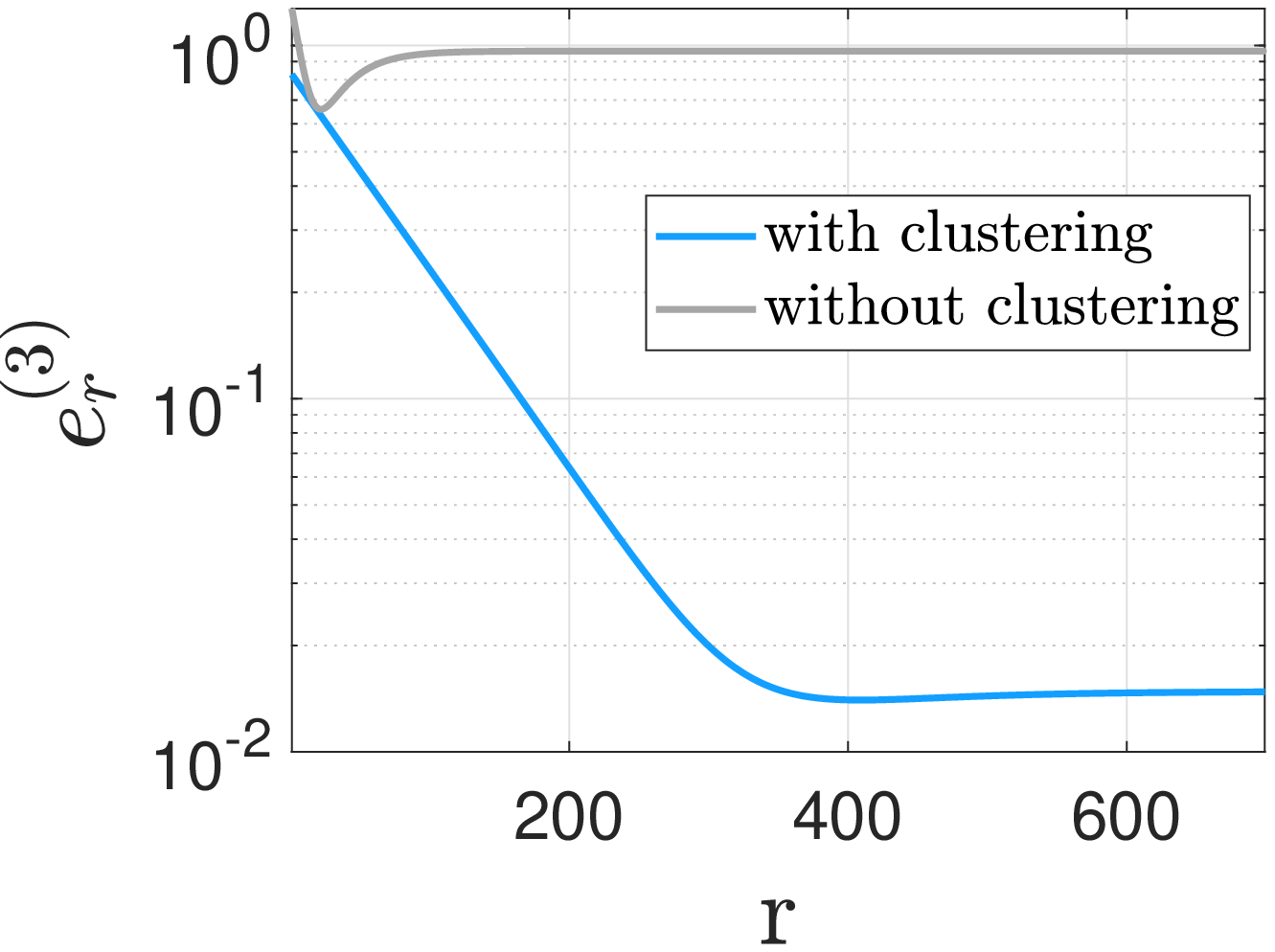}\label{fig:with_without_cluster_3} \vspace{-1em}
   \includegraphics[width=0.3\textwidth]{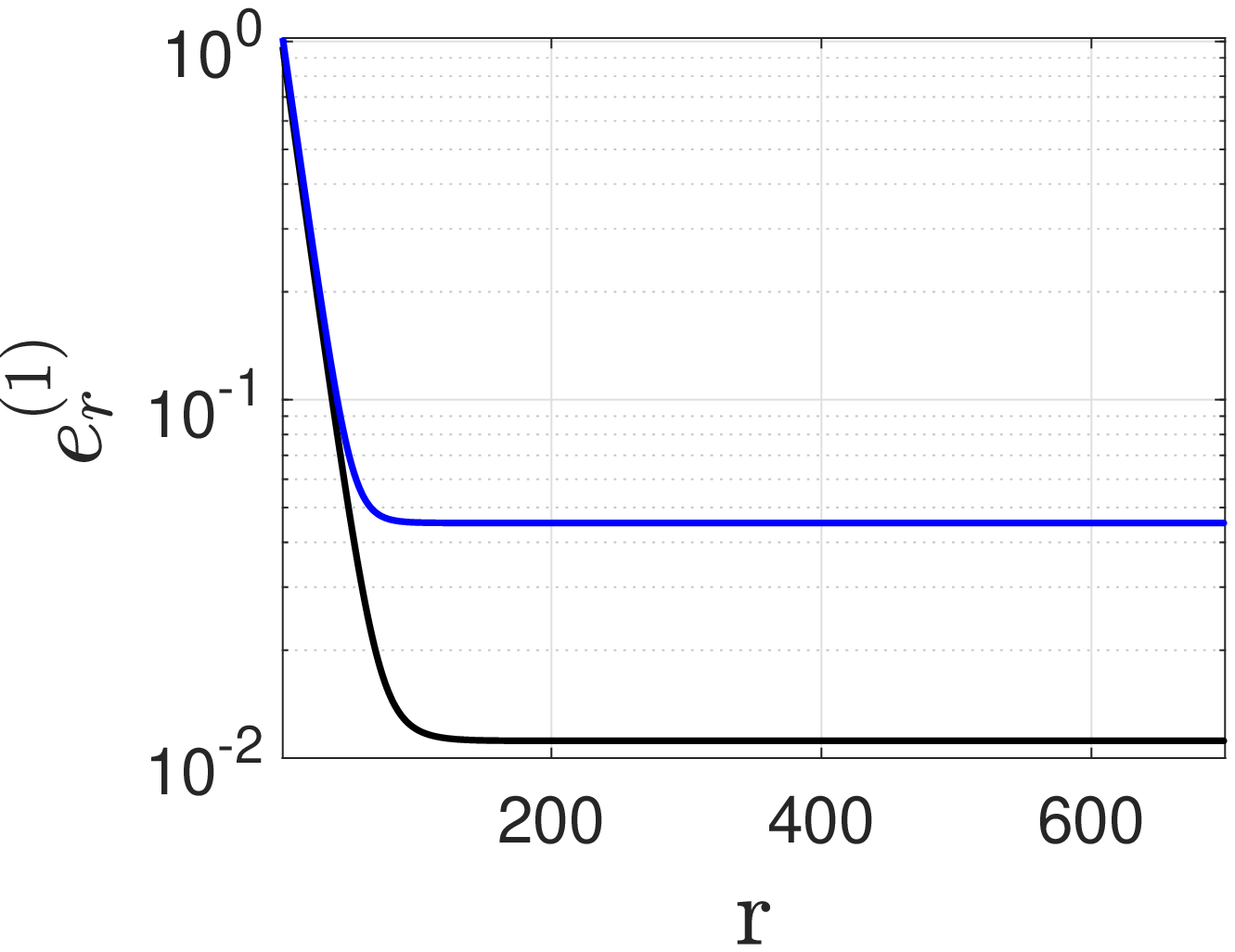}\label{fig:single_multiple_cluster_1}
    \includegraphics[width=0.3\textwidth]{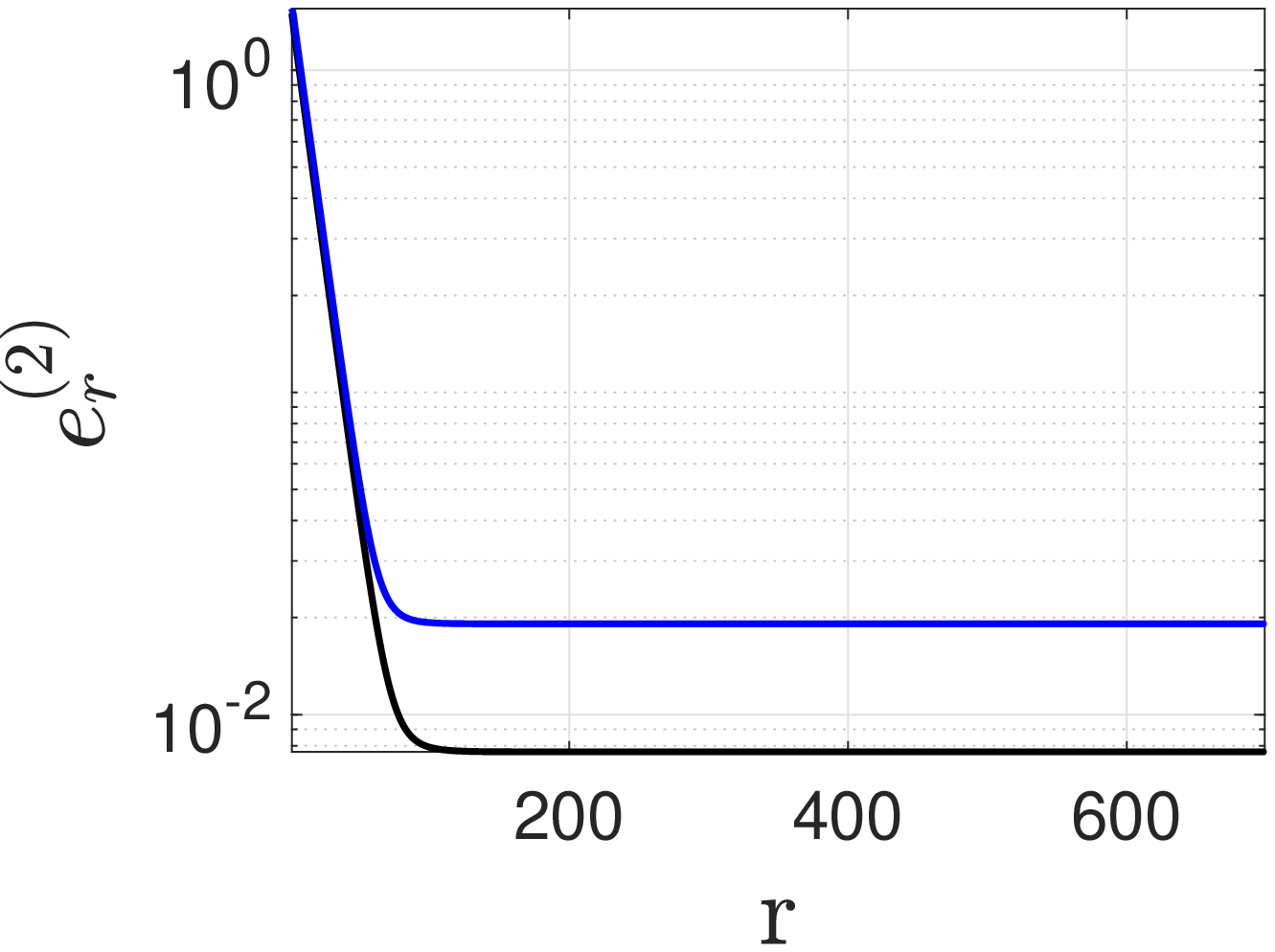}\label{fig:single_multiple_cluster_2}
   \includegraphics[width=0.3\textwidth]{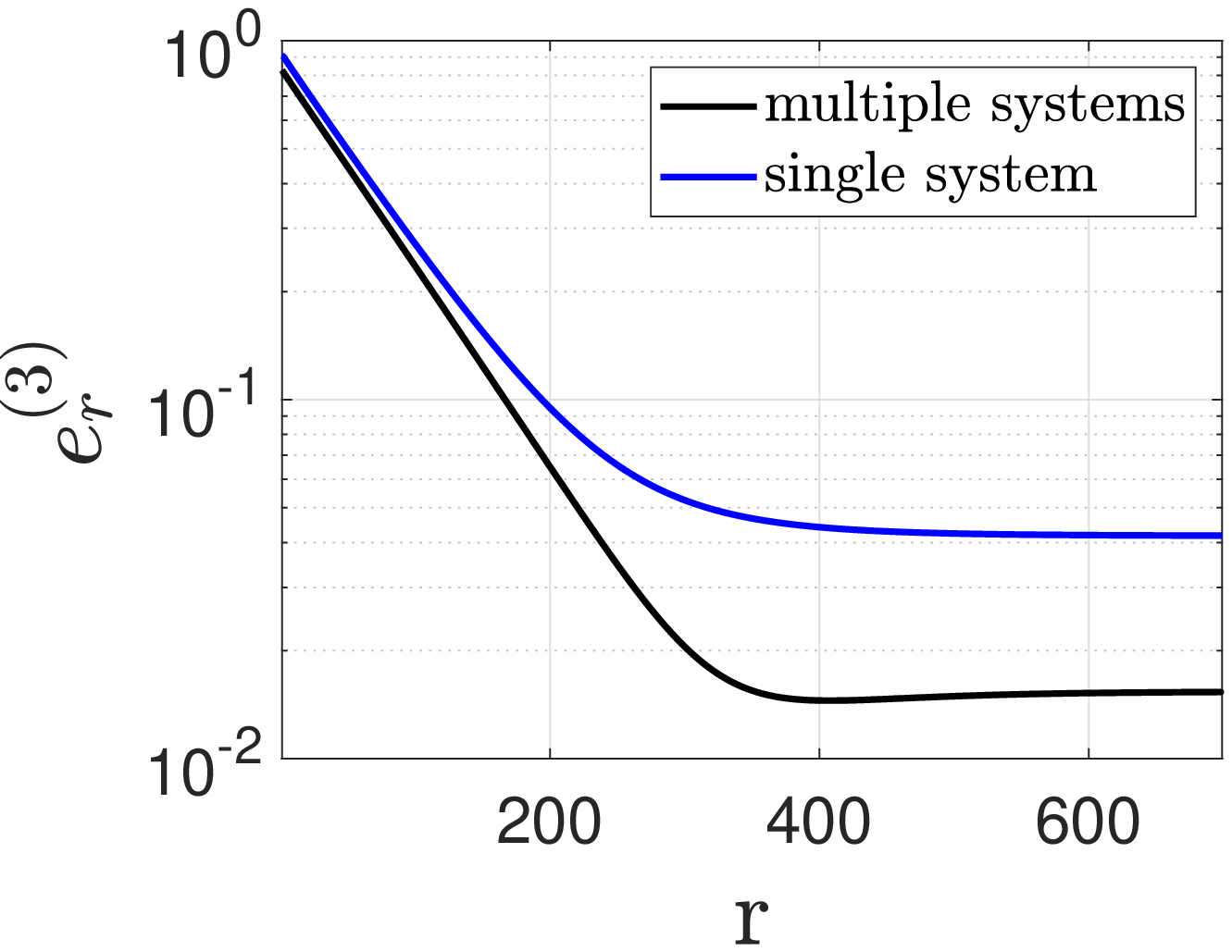}\label{fig:single_multiple_cluster_3}
   \vspace{-0.2em}
    \caption{Estimation error as a function of iteration count. The plot on the top considers Algorithm \ref{algorithm:clustered_fedsysid} with and without clustering, whereas the bottom plot consider the single and multiple agents settings.}
    \label{fig:convergence}

\end{figure*}

Figure \ref{fig:convergence} depicts the estimation error $e^{(j)}_r$ as a function of the number of iterations $r$ for all the three considered clusters. The top plots compare the performance of Algorithm \ref{algorithm:clustered_fedsysid} with and without the clustering procedure (i.e., line 5 of Algorithm \ref{algorithm:clustered_fedsysid}). These plots reveals that the estimation error decreases significantly when systems with the same model are clustered and cooperate to estimate their dynamics. Conversely, in the absence of  clustering, the significant heterogeneity level across the systems leads to a poor common estimation, resulting in a large estimation error and unpersonalized solutions. This confirms our theoretical results, showing that the misclassification rate in \eqref{eq:epsilon_corollary} outperforms the heterogeneity constant of \citep{xin2022identifying,xin2023learning,wang2022fedsysid}, when dealing with heterogeneous settings.

The bottom plots of Figure \ref{fig:convergence} demonstrates the benefits of collaboration among systems to learn their dynamics. This shows that the estimation error is considerably reduced when multiple systems within the same cluster (i.e., $|\mathcal{C}_1|=10$, $|\mathcal{C}_2|=24$, and $|\mathcal{C}_3|=16$) leverage the data from each other to identify their dynamics, compared to the case where a single system estimate its dynamics by using its own observations. This also confirms our theoretical results, where the statistical error in \eqref{eq:epsilon_corollary} scales down with the number of systems in the cluster, thus highlighting the benefit of collaboration in improving estimation accuracy in a multi-system setting. 

\begin{figure}[H]
    \centering
   \includegraphics[width=0.5\textwidth]{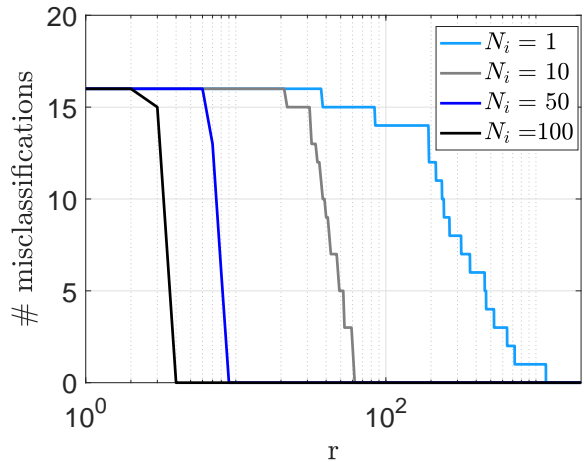}
    \caption{Number of misclassification as a function  of iteration count.}
    \label{fig:misclassification}
\end{figure}

Figure \ref{fig:misclassification} presents the misclassifications of Algorithm \ref{algorithm:clustered_fedsysid} as a function of  iterations $r$. It depicts the number of systems whose cluster identity is incorrectly estimated. The figure illustrates the effect of the number of observed trajectories on the misclassification rate. As anticipated and consistent with our theoretical results, an increase in the number of trajectories leads to a considerable reduction in the number of iterations needed to correctly classify all the systems into their respective clusters.  

\section{Conclusions and Future Work} \label{sec:conclusions}

We presented an approach to address the system identification problem through the use of clustering. Our method involves partitioning different systems that observe multiple trajectories into disjoint clusters based on the similarity of their dynamics. This approach enjoys an improved convergence rate that scales inversely with the number of systems in the cluster, along with an additive misclassification rate that has been shown to be negligible under mild assumptions. Our approach enables systems within the same cluster to learn their dynamics more efficiently. 
Future work will involve extending the proposed formulation to online system identification and proposing an adaptive clustering approach that eliminates the necessity for the warm initialization and well-separated clusters assumptions.

\bibliography{bibliography.bib}
\bibliographystyle{abbrvnat}

\newpage

\section{Appendix}

\subsection{Proof of Lemma \ref{lemma:prob_misclassification}} \label{proof:probability_of_misclassification}

Without loss of generality, we can analyze only the first cluster $\mathcal{M}_i^{1, j}$ for some $j \neq 1$. By definition, we have
$$
\mathcal{M}_i^{1, j}=\left\{ \| X^{(i)} - \widehat{\Theta}_jZ^{(i)}\|_F^2 \leq \| X^{(i)} - \widehat{\Theta}_1Z^{(i)}\|_F^2 \right\}
$$
where the batch matrices $X^{(i)}, Z^{(i)}$ and $ W^{(i)}$ are related according to $X^{(i)}={\Theta}_1 Z^{(i)}+W^{(i)}$. Note that $z^{(i)}_{l,t} \stackrel{\text{i.i.d.}}{\sim} \mathcal{N}\left(0, \Sigma^{(i)}_t\right)$ and $w^{(i)}_{l,t} \stackrel{\text{i.i.d.}}{\sim} \mathcal{N}\left(0, \sigma_{w,i}^{2} I_{n_x}\right)$ are independent across trajectories (i.e., the columns of $Z^{(i)}$ and $W^{(i)}$ are independent). Thus, we can write

\begin{align*}
\mathbb{P}\left\{\mathcal{M}_i^{1, j}\right\}&=\mathbb{P}\left\{\left\|({\Theta}_1 - \widehat{\Theta}_1)Z^{(i)} + W^{(i)}\right\|_F^2 \geq\left\|({\Theta}_1 - \widehat{\Theta}_j)Z^{(i)} + W^{(i)}\right\|_F^2\right\}\\
&= \mathbb{P}\left\{\sum_{t=0}^{T-1}\sum_{l=1}^{{N_i}} m^{(i),\top}_{l,t}m^{(i)}_{l,t} \geq \sum_{t=0}^{T-1}\sum_{l=1}^{{N_i}} n^{(i),\top}_{l,t}n^{(i)}_{l,t} \right\},
\end{align*}

where $m^{(i)}_{l,t} = ({\Theta}_1 - \widehat{\Theta}_1)z^{(i)}_{l,t} + w^{(i)}_{l,t}  \sim \mathcal{N}\left(0, \bar{\Sigma}^{(i)}_t\right) $  , $n^{(i)}_{l,t} = ({\Theta}_1 - \widehat{\Theta}_j)z^{(i)}_{l,t} + w^{(i)}_{l,t} \sim \mathcal{N}\left(0, \tilde{\Sigma}^{(i)}_t\right)$, with
\begin{align*}
    \bar{\Sigma}^{(i)}_t &= ({\Theta}_1 - \widehat{\Theta}_1)\Sigma^{(i)}_t ({\Theta}_1 - \widehat{\Theta}_1)^\top + \sigma^2_{w,i}I_{n_x},\;\
    \Tilde{\Sigma}^{(i)}_t =({\Theta}_1 - \widehat{\Theta}_j)\Sigma^{(i)}_t ({\Theta}_1 - \widehat{\Theta}_j)^\top + \sigma^2_{w,i}I_{n_x}.
\end{align*}
 Therefore, we obtain
\begin{align*}
\mathbb{P}\left\{\mathcal{M}_i^{1, j}\right\}
&= \mathbb{P}\left\{\sum_{t=0}^{T-1}\sum_{l=1}^{{N_i}} v^{(i),\top}_{l,t}\bar{\Sigma}^{(i)}_tv^{(i)}_{l,t} \geq \sum_{t=0}^{T-1}\sum_{l=1}^{{N_i}} u^{(i),\top}_{l,t}\Tilde{\Sigma}^{(i)}_t  u^{(i)}_{l,t} \right\},
\end{align*}
with $m^{(i)}_{l,t}=(\bar{\Sigma}^{(i)}_t)^{\frac{1}{2}} v^{(i)}_{l,t}$ and $n^{(i)}_{l,t}=(\bar{\Sigma}^{(i)}_t)^{\frac{1}{2}} u^{(i)}_{l,t}$ for some standard normal random vectors
$v^{(i)}_{l,t}$, $u^{(i)}_{l,t}$ $\sim$  $\mathcal{N}\left(0, I_{n_x}\right)$. Then, the above expression can be rewritten as follows
\begin{align*}
\mathbb{P}\left\{\mathcal{M}_i^{1, j}\right\}
&= \mathbb{P}\left\{\sum_{t=0}^{T-1}\sum_{l=1}^{{{N_i}}} v^{(i),\top}_{l,t}\bar{\Sigma}^{(i)}_tv^{(i)}_{l,t} \geq \sum_{t=0}^{T-1}\sum_{l=1}^{{{N_i}}} \|\Tilde{\Sigma}^{(i)}_t\|  u^{(i),\top}_{l,t} u^{(i)}_{l,t} \right\}\\
&=\mathbb{P}\left\{\sum_{t=0}^{T-1}\sum_{l=1}^{{{N_i}}} v^{(i),\top}_{l,t}\bar{\Sigma}^{(i)}_tv^{(i)}_{l,t} \geq \sum_{t=0}^{T-1}\sum_{l=1}^{{{N_i}}} c^{(i)}_t  u^{(i),\top}_{l,t} u^{(i)}_{l,t} \right\}
\end{align*}
with $c^{(i)}_t=\|{\Theta}_1 - \widehat{\Theta}_j\|^2 \|\Sigma^{(i)}_t\| + \sigma^2_{w,i}\sqrt{n_x}$, which implies
\begin{align*}
\mathbb{P}\left\{\mathcal{M}_i^{1, j}\right\}
&= \mathbb{P}\left\{\sum_{t=0}^{T-1}\sum_{l=1}^{{{N_i}}} v^{(i),\top}_{l,t}\bar{\Sigma}^{(i)}_tv^{(i)}_{l,t} \geq \sum_{t=0}^{T-1}\sum_{l=1}^{{{N_i}}} c^{(i)}_t  u^{(i),\top}_{l,t} u^{(i)}_{l,t} \right\}\\
& \leq \mathbb{P}\left\{\sum_{t=0}^{T-1}\sum_{l=1}^{{{N_i}}} c^{(i)}_t  u^{(i),\top}_{l,t} u^{(i)}_{l,t} \leq \bar{t}\right\} +\mathbb{P}\left\{\sum_{t=0}^{T-1}\sum_{l=1}^{{{N_i}}} v^{(i),\top}_{l,t}\bar{\Sigma}^{(i)}_tv^{(i)}_{l,t} > \bar{t} \right\},
\end{align*}
for any $\bar{t}\geq 0$. Therefore , by using $v^{(i),\top}_{l,t}\bar{\Sigma}^{(i)}_tv^{(i)}_{l,t} \leq d^{(i)}_t v^{(i),\top}_{l,t} v^{(i)}_{l,t}$ with $d^{(i)}_t=\|{\Theta}_1 - \widehat{\Theta}_1\|^2 \|\Sigma^{(i)}_t\| + \sigma^2_{w,i}\sqrt{n_x}$ we obtain 
\begin{align*}
\mathbb{P}\left\{\mathcal{M}_i^{1, j}\right\}
\leq \mathbb{P}\left\{\sum_{t=0}^{T-1} c^{(i)}_t  V^{(i)}_t \leq \bar{t}\right\}+\mathbb{P}\left\{\sum_{t=0}^{T-1}d^{(i)}_t V^{(i)}_t > \bar{t} \right\},
\end{align*}
where $V^{(i)}_t$ are standard Chi-squared distributions with $N_i n_x$ degrees of freedom, for all $t \in \{0,1,\ldots, T-1\}$. Moreover, by using Definition \ref{def:minimum_separation} and Assumption \ref{assumption:warm_initialization}, 
\begin{align*}
\mathbb{P}\left\{\mathcal{M}_i^{1, j}\right\}
\leq \mathbb{P}\left\{\sum_{t=0}^{T-1} f^{(i)}_t  V^{(i)}_t \leq \bar{t}\right\}+\mathbb{P}\left\{\sum_{t=0}^{T-1}g^{(i)}_t V^{(i)}_t > \bar{t} \right\},
\end{align*}
with $f^{(i)}_t = (\frac{1}{2} + \alpha)^2 \Delta_{\min}^2 \|\Sigma^{(i)}_t\| + \sigma^2_{w,i}\sqrt{n_x}$ and $g^{(i)}_t = (\frac{1}{2} - \alpha)^2 \Delta_{\min}^2 \|\Sigma^{(i)}_t\| + \sigma^2_{w,i}\sqrt{n_x}$, since $c^{(i)}_t = \|{\Theta}_1 - \widehat{\Theta}_j\|^2 \|\Sigma^{(i)}_t\| + \sigma^2_{w,i}\sqrt{n_x}\geq (\frac{1}{2} + \alpha)^2 \Delta_{\min}^2 \|\Sigma^{(i)}_t\| + \sigma^2_{w,i}\sqrt{n_x},$ with $ \|{\Theta}_j - \widehat{\Theta}_1\| \geq \|{\Theta}_j - {\Theta}_1\| - \|\widehat{\Theta}_j - {\Theta}_j \| = (\frac{1}{2}+\alpha)\Delta_{\min}$ and $d^{(i)}_t = \|{\Theta}_1 - \widehat{\Theta}_1\|^2 \|\Sigma^{(i)}_t\| + \sigma^2_{w,i}\sqrt{n_x} 
\leq (\frac{1}{2} - \alpha)^2 \Delta_{\min}^2 \|\Sigma^{(i)}_t\| + \sigma^2_{w,i}\sqrt{n_x}$, where  $\|{\Theta}_1 - \widehat{\Theta}_1\| \leq (\frac{1}{2}- \alpha)\Delta_{\min}$ according to Assumption \ref{assumption:warm_initialization}. Therefore, to characterize the above tail bounds, we can exploit well-established concentration inequalities as detailed in \citep{boucheron2013concentration, vershynin2010introduction}. To this end, we can use union bound to write
\begin{align*}
\mathbb{P}\left\{\mathcal{M}_i^{1, j}\right\}
\leq \sum_{t=0}^{T-1}\mathbb{P}\left\{ f^{(i)}_t  V^{(i)}_t \leq \bar{t}\right\}+\mathbb{P}\left\{g^{(i)}_t V^{(i)}_t > \bar{t} \right\},
\end{align*}
where $\mathbb{P}\left\{ f^{(i)}_t  V^{(i)}_t \leq \bar{t}\right\}$ can be rewritten as follows
\begin{align*}
\mathbb{P}\left\{  f^{(i)}_t V^{(i)}_t \leq \bar{t}\right\} = \mathbb{P}\left\{ V^{(i)}_t \leq \frac{4\bar{t}}{\sigma^2_{w,i}\sqrt{n_x}\left((1+2\alpha)^2\rho^{(i)} \frac{\| \Sigma^{(i)}_t\|}{\sqrt{n_x}} + 4\right)}\right\},
\end{align*}
thus, by choosing $\Bar{t}=N_in_x\left((\frac{1}{4}+\alpha^2)\Delta^2_{\min}\|\Sigma^{(i)}_t\| + \sigma^2_{w,i}\sqrt{n_x}\right)$ we obtain
\begin{align*}
\mathbb{P}\left\{ f^{(i)}_t V^{(i)}_t\leq \bar{t}\right\} = \mathbb{P}\left\{ \frac{V^{(i)}_t}{N_i n_x} - 1  \leq \frac{-4\alpha \|\Sigma^{(i)}_t\|}{(1+2\alpha)^2\rho^{(i)} \| \Sigma^{(i)}_t\|+ 4\sqrt{n_x}}\right\},
\end{align*}
\noindent as per the concentration of standard Chi-squared distributions in \citep{wainwright2019high}, it is established that there exist universal constants $c_1$ and $c_2$, such that
\begin{align} \label{eq:prob_fV}
\mathbb{P}\left\{ f^{(i)}_t V^{(i)}_t \leq \bar{t}\right\} 
& \leq c_1 \exp\left(-c_2 N_in_x \left(\frac{\alpha \rho^{(i)} \|\Sigma^{(i)}_t\|}{\rho^{(i)}\|\Sigma^{(i)}_t\| + \sqrt{n_x} }\right)^2\right).
\end{align}

Similarly, $\mathbb{P}\left\{g^{(i)}_t V^{(i)}_t > \bar{t} \right\}$ can be rewritten as follows
\begin{align*}
\mathbb{P}\left\{ g^{(i)}_t V^{(i)}_t\leq \bar{t}\right\} = \mathbb{P}\left\{ \frac{V^{(i)}_t}{N_i n_x} - 1  \leq \frac{4\alpha \|\Sigma^{(i)}_t\|}{(1-2\alpha)^2\rho^{(i)} \| \Sigma^{(i)}_t\|+ 4\sqrt{n_x}}\right\},
\end{align*}
and by the concentration of Chi-squared distribution
\begin{align}\label{eq:prob_gV}
\mathbb{P}\left\{ g^{(i)}_t V^{(i)}_t \leq \bar{t}\right\} &\leq c_3 \exp\left(-c_4 N_in_x \left(\frac{\alpha \rho^{(i)} \|\Sigma^{(i)}_t\|}{\rho^{(i)}\|\Sigma^{(i)}_t\| + \sqrt{n_x} }\right)^2\right),
\end{align}
where the proof is completed by combining \eqref{eq:prob_fV} and \eqref{eq:prob_gV} to obtain
 \begin{align*}
\mathbb{P}\left\{\mathcal{M}_i^{1, j}\right\}
\leq c_1 \sum_{t=0}^{T-1} \exp\left(-c_2 N_in_x \left(\frac{\alpha \rho^{(i)} \|\Sigma^{(i)}_t\|}{\rho^{(i)}\|\Sigma^{(i)}_t\| + \sqrt{n_x} }\right)^2\right).
\end{align*}

\subsection{Proof of Theorem \ref{theorem:convergence}} \label{proof:convergence}

 Without loss of generality, we analyze only the first cluster. Recall that the model is updated as follows:
\begin{align}\label{eq:model_update_rule}
\widehat{\Theta}^+_1&=\frac{1}{|\widehat{\mathcal{C}}_1|}\sum_{i \in \widehat{\mathcal{C}}_1 }\tilde{\Theta}_i
=\frac{1}{|\widehat{\mathcal{C}}_1|} \sum_{i \in \widehat{\mathcal{C}}_1 \cap \mathcal{S}_1 }\tilde{\Theta}_i + \frac{1}{|\widehat{\mathcal{C}}_1|} \sum_{i \in \widehat{\mathcal{C}}_1 \cap \overline{\mathcal{S}_1} }\tilde{\Theta}_i
\end{align}
with $\Tilde{\Theta}_i = \widehat{\Theta}_{1} + 2\eta_1(X^{(i)} - \widehat{\Theta}_{1}Z^{(i)})Z^{(i),\top}$. Here $\widehat{\mathcal{C}}_1 \cap \mathcal{C}_1$ corresponds to the set of systems correctly classified to the first cluster and $\widehat{\mathcal{C}}_1 \cap \overline{\mathcal{C}_1}$ represents the set of systems that are misclassified to the first cluster, with $\overline{\mathcal{C}_1}$ denoting the complement of $\mathcal{C}_1$. The above expression can be rewritten as follows
\begin{align*}
\widehat{\Theta}^+_1&=\widehat{\Theta}_1 + \frac{2\eta_1}{|\widehat{\mathcal{C}}_1|} \sum_{i \in \widehat{\mathcal{C}}_1 \cap \mathcal{C}_1} (X^{(i)} - \widehat{\Theta}_{1}Z^{(i)})Z^{(i),\top}
+ \frac{2\eta_1}{|\widehat{\mathcal{C}}_1|} \sum_{i \in \widehat{\mathcal{C}}_1 \cap \overline{\mathcal{C}_1}} (X^{(i)} - \widehat{\Theta}_{1}Z^{(i)})Z^{(i),\top},
\end{align*}
where $X^{(i)} = \Theta_1 Z^{(i)} + W^{(i)}$ for $i \in \widehat{\mathcal{C}}_1 \cap \mathcal{C}_1$, and $X^{(i)} = \Theta_j Z^{(i)} + W^{(i)}$ for $i \in \widehat{\mathcal{C}}_1 \cap \overline{\mathcal{C}_1}$, with $j \neq 1 \in [K]$. Therefore, by manipulating the above expression, we have
\begin{align*}
\widehat{\Theta}^+_1-\Theta_1 &= (\widehat{\Theta}_1-\Theta_1)\left(I  - \frac{2\eta_1}{|\widehat{\mathcal{C}}_1|} \sum_{i \in \widehat{\mathcal{C}}_1} Z^{(i)}Z^{(i),\top} \right)
+ \frac{2 \eta_1}{|\widehat{\mathcal{C}}_1|} \sum_{i \in \widehat{\mathcal{C}}_1} W^{(i)}Z^{(i),\top}\\
&+ (\Theta_j - \Theta_1) \frac{2\eta_1}{|\widehat{\mathcal{C}}_1|} |\widehat{\mathcal{C}}_1 \cap \overline{\mathcal{C}_1}| \sum_{i \in \widehat{\mathcal{C}}_1 \cap \overline{\mathcal{C}_1}} Z^{(i)}Z^{(i),\top},
\end{align*}   
and thus, we obtain
\begin{align*}
\|\widehat{\Theta}^+_1 - {\Theta}_1\| \leq & \| \mathcal{H}_1\| + \|\mathcal{H}_2\|,
\end{align*}
with,
\begin{align*}
    \|\mathcal{H}_1\| & = \|\widehat{\Theta}_1-\Theta_1\| \left\|I  - \frac{2\eta_1}{|\widehat{\mathcal{C}}_1|}  ZZ^{\top}\right\|+ \frac{2 \eta_1}{|\widehat{\mathcal{C}}_1|} \sum_{i \in \widehat{\mathcal{C}}_1}\| WZ^{\top}\|,
\end{align*}
\begin{align*}
    \|\mathcal{H}_2\| & = \|\Theta_j - \Theta_1\| \frac{2\eta_1}{|\widehat{\mathcal{C}}_1|} |\widehat{\mathcal{C}}_1 \cap \overline{\mathcal{C}_1}| \|\bar{Z}\bar{Z}^{\top}\|.
\end{align*}

We now concatenate the batch matrices $Z^{(i)}, W^{(i)}$ of the systems classified to the first cluster in $Z \in \mathbb{R}^{(n_x+n_u)\times N_i T |\widehat{\mathcal{C}}_1|}$ and $W \in \mathbb{R}^{n_x\times N_i T |\widehat{\mathcal{C}}_1|}$, and similarly the batch matrices $Z^{(i)}$ of the systems incorrectly classified to the first cluster are concatenated in $\Bar{Z} \in \mathbb{R}^{(n_x+n_u)\times N_i T |\widehat{\mathcal{C}}_1 \cap \overline{\mathcal{C}_1}|}$. We proceed with our analysis by controlling both terms separately. To upper bound the first term, we introduce the following propositions.

\begin{prop} \citep[Proposition 8]{wang2022fedsysid}
For any fixed $0<\delta<1$, let $N_{i} \geq (4 n_x+$
$2 n_u) \log \frac{T|\widehat{\mathcal{C}}_1|}{\delta}$. It holds, with probability at least $1-\delta$, that
\begin{align} \label{eq:prop_WZ}
 \left\|W Z^{\top}\right\| \leq  4 \sigma_{w,i} \sqrt{N_{i}(2 n_x+n_u) \log \frac{9|\widehat{\mathcal{C}}_1| T}{\delta}} \sum_{t=0}^{T-1} \left\|(\Sigma_{t}^{(i)})^{\frac{1}{2}}\right\| . 
\end{align}
\end{prop}

\begin{prop} (Adapted from \citep[Proposition 6 ]{wang2022fedsysid}) For any fixed $0<\delta<1$, let $N_{i} \geq 8(n_x+ n_u)+16 \log \frac{2|\widehat{\mathcal{C}}_1|T}{\delta}$. It holds, with probability at least $1-\delta$, that
\begin{align}\label{eq:prop_ZZ}
    ZZ^\top \succeq  \frac{1}{4}\sum_{i \in \widehat{\mathcal{C}}_1 } N_i \sum_{t=0}^{T-1} {\Sigma}_{t}^{(i)},
\end{align}
\begin{align}\label{eq:prop_ZZbar}
\|\bar{Z}\bar{Z}^\top\| \leq \frac{9}{4}\sum_{i \in \widehat{\mathcal{C}}_1 \cap \overline{\mathcal{C}_1}}  \sum_{t=0}^{T-1}N_i \left\|\Sigma_{t}^{(i)}\right\|.
\end{align}
\begin{proof}
Expression~\eqref{eq:prop_ZZ} follows direct from Proposition 6 in \citep{wang2022fedsysid}. For expression~\eqref{eq:prop_ZZbar}, we can first write
\begin{align*}
 \|\bar{Z}\bar{Z}^\top\| &= \left\|  \sum_{i \in \widehat{\mathcal{C}}_1 \cap \overline{\mathcal{C}_1} } \sum_{l=1}^{N_i} \sum_{t=0}^{T-1} z^{(i)}_{l,t}z^{(i),\top}_{l,t}     \right\|\\
 &\leq  \sum_{i \in \widehat{\mathcal{C}}_1 \cap \overline{\mathcal{C}_1} } \left\| \sum_{l=1}^{N_i} \sum_{t=0}^{T-1} z^{(i)}_{l,t}z^{(i),\top}_{l,t}     \right\|
\end{align*}
where ${\chi}_{l,t}^{(i)}=(\Sigma_{t}^{(i)})^{-\frac{1}{2}} z_{l,t}^{(i)}$ for any fixed $l,t$, and $i$, where ${\chi}_{l,t}^{i} \stackrel{\text{i.i.d.}}{\sim} \mathcal{N}\left(0, I_{n_x+n_u}\right)$, for all $l \in$ $\left\{1,2, \ldots, N_{i}\right\}$, we obtain 
\begin{align*}
 \|\bar{Z}\bar{Z}^\top\| \leq  \sum_{i \in \widehat{\mathcal{C}}_1 \cap \overline{\mathcal{C}_1} } \sum_{t=0}^{T-1} \|\Sigma^{(i)}_t\| \left\| \sum_{l=1}^{N_i}\chi^{(i)}_{l,t}\chi^{(i),\top}_{l,t} \right\|,
\end{align*}
thus, by using Proposition 6 of \citep{wang2022fedsysid}, with probability $1-\frac{\delta}{T}$, we have  
\begin{align*}
    \left\| \sum_{l=1}^{N_i}\chi^{(i)}_{l,t}\chi^{(i),\top}_{l,t}     \right\| \leq \frac{9N_i}{4},
\end{align*}
which implies
\begin{align*}
 \|\bar{Z}\bar{Z}^\top\| \leq \frac{9}{4} \sum_{i \in \widehat{\mathcal{C}}_1 \cap \overline{\mathcal{C}_1} }\sum_{t=0}^{T-1} N_i  \|\Sigma^{(i)}_t\|.
\end{align*}
\end{proof}
\end{prop}

Therefore, with probability $1-2\delta$, we have
\begin{align*}
\|\mathcal{H}_1\| &\leq \|\widehat{\Theta}_1-\Theta_1\| \left\|I  - \frac{\eta_1}{2|\widehat{\mathcal{C}}_1|}\sum_{i \in \widehat{\mathcal{C}}_1 } N_i \sum_{t=0}^{T-1} {\Sigma}_{t}^{(i)}      \right\|
+ \frac{2 \eta_1}{|\widehat{\mathcal{C}}_1|} \sum_{i \in \widehat{\mathcal{C}}_1}\| WZ^{\top}\|,\\
&=\|\widehat{\Theta}_1-\Theta_1\| \left(1  - \frac{\eta_1}{2|\widehat{\mathcal{C}}_1|}\lambda_{\min}\left(\sum_{i \in \widehat{\mathcal{C}}_1 } N_i \sum_{t=0}^{T-1} {\Sigma}_{t}^{(i)}\right)    \right)
+ \frac{2 \eta_1}{|\widehat{\mathcal{C}}_1|} \sum_{i \in \widehat{\mathcal{C}}_1}\| WZ^{\top}\|.
\end{align*} 

Hence, by selecting $\eta_1 =\frac{|\widehat{\mathcal{C}}_1|}{\lambda_{\min}\left(\sum_{i \in \widehat{\mathcal{C}}_1 } N_i \sum_{t=0}^{T-1} {\Sigma}_{t}^{(i)}\right)}$, we obtain
\begin{align}\label{eq:bound_H1}
\|\mathcal{H}_1\| &\leq  \frac{1}{2} \|\widehat{\Theta}_1-\Theta_1\|
+ \frac{8\sum_{i \in \widehat{\mathcal{C}}_1} \sigma_{w,i} \sqrt{N_{i}(2 n_x+n_u) \log \frac{9|\widehat{\mathcal{C}}_1|T}{\delta}} \sum_{t=0}^{T-1} \left\|(\Sigma_{t}^{(i)})^{\frac{1}{2}}\right\|}{\lambda_{\min }\left(N_{i} \sum_{t=0}^{T-1} {\Sigma}^{(i)}_t\right)}\notag\\
& \leq \frac{1}{2} \|\widehat{\Theta}_1-\Theta_1\|
+\frac{8 \sqrt{(2 n_x+n_u) \log \frac{9|\widehat{\mathcal{C}}_1| T}{\delta}}\sqrt{\sum_{i \in \mathcal{S}_1} \sigma_{w,i}^2 \left(\sum_{t=0}^{T-1}\left\|({\Sigma}_{t}^{(i)})^{\frac{1}{2}}\right\|\right)^2}}{\sqrt{\sum_{i \in \widehat{\mathcal{C}}_1} N_i}\times \min_{i \in \widehat{\mathcal{C}}_1}\lambda_{\min }\left(\sum_{t=0}^{T-1}{\Sigma}^{(i)}_t\right)}\notag \\
&= \frac{1}{2}  \|\widehat{\Theta}_1-\Theta_1\|+\bar{c}_0 \times \frac{1}{\sqrt{\sum_{i \in \widehat{\mathcal{C}}_1} N_i}},
\end{align}
with $N_{i} \geq$ $\max \{8(n_x+n_u)+16 \log \frac{2|\widehat{\mathcal{C}}_1| T}{\delta}$, $ (4 n_x+2 n_u) \log \frac{|\widehat{\mathcal{C}}_1| T}{\delta} \}$, for all $i\in \widehat{\mathcal{C}}_1$. To control the second term $\|\mathcal{H}_2\|$, we first use the Definition \ref{def:minimum_separation} to write
\begin{align*}
\|\mathcal{H}_2\| \leq \Delta_{\max} |\widehat{\mathcal{C}}_1 \cap \overline{\mathcal{C}_1}| \frac{9\sum_{i \in \widehat{\mathcal{C}}_1 \cap \overline{\mathcal{C}_1} } N_i \sum_{t=0}^{T-1} \|\Sigma^{(i)}_t\|}{2\lambda_{\min}\left(\sum_{i \in \widehat{\mathcal{C}}_1 } N_i \sum_{t=0}^{T-1} {\Sigma}_{t}^{(i)}\right)} ,
\end{align*}
which implies
\begin{align*}
\|\mathcal{H}_2\| \leq c_5\Delta_{\max} |\widehat{\mathcal{C}}_1 \cap \overline{\mathcal{C}_1}|  ,
\end{align*}
by using Jensen and Cauchy-Schwartz inequalities in the denominator and numerator, respectively, where we define $c_5=\frac{9\sum_{i \in \widehat{\mathcal{C}}_1 \cap \overline{\mathcal{C}_1} }  \sum_{t=0}^{T-1} \|\Sigma^{(i)}_t\|}{2\min_{i \in \widehat{\mathcal{C}}_1 }\left( \sum_{t=0}^{T-1} {\Sigma}_{t}^{(i)}\right)}$. Therefore, we proceed with our analysis to control $|\widehat{\mathcal{C}}_1 \cap \overline{\mathcal{C}_1}|$. To do so, we use Lemma \ref{lemma:prob_misclassification} and obtain
\begin{align*}
\mathbb{E}\left[|\widehat{\mathcal{C}}_1 \cap \overline{\mathcal{C}_1}|\right] &\leq c_6\sum_{i \in [M]} \sum_{t=0}^{T-1} \exp\left(-c_7 N_in_x \left(\frac{\alpha \rho^{(i)} \|\Sigma^{(i)}_t\|}{\rho^{(i)}\|\Sigma^{(i)}_t\| + \sqrt{n_x} }\right)^2\right),\\
\end{align*}
which yields
\begin{align*}
&\mathbb{P}\left\{|\widehat{\mathcal{C}}_1 \cap \overline{\mathcal{C}_1}| \leq c_6 \sum_{i \in [M]} \sum_{t=0}^{T-1} \exp\left(-\frac{c_7}{2} N_in_x \left(\frac{\alpha \rho^{(i)} \|\Sigma^{(i)}_t\|}{\rho^{(i)}\|\Sigma^{(i)}_t\| + \sqrt{n_x} }\right)^2\right) \right\}\\
&\geq 1 - \sum_{i \in [M]} \sum_{t=0}^{T-1} \exp\left(-\frac{c_7}{2} N_in_x \left(\frac{\alpha \rho^{(i)} \|\Sigma^{(i)}_t\|}{\rho^{(i)}\|\Sigma^{(i)}_t\| + \sqrt{n_x} }\right)^2\right)  \geq 1 -\delta, 
\end{align*}    
by using Markov's inequality and  Assumption \ref{assumption:Ni_nx} with $N_in_x \geq c\left(\frac{\rho^{(i)}\|\Sigma^{(i)}_t\| + \sqrt{n_x}}{\alpha \rho^{(i)} \|\Sigma^{(i)}_t\| }\right)^2
\log(\frac{MT}{\delta})$, for some large enough constant $c$ such that $\frac{1}{c} < c_7$, with $0<\delta<1$ for all $i \in [M]$. Thus, we obtain 
\begin{align}\label{eq:bound_H2}
\|\mathcal{H}_2\| \leq \bar{c}_1\Delta_{\max} \sum_{i \in [M] }\sum_{t=0}^{T-1} \exp\left(-\bar{c}_2 N_in_x \left(\frac{\alpha \rho^{(i)} \|\Sigma^{(i)}_t\|}{\rho^{(i)}\|\Sigma^{(i)}_t\| + \sqrt{n_x} }\right)^2\right) ,
\end{align}
with probability at least $1-\delta$. The proof is completed by combining \eqref{eq:bound_H1} and \eqref{eq:bound_H2}.

\subsection{Proof of Corollary \ref{corollary:convergence}} \label{proof:corollary_convergence}

We first recall that at iteration $r$ we posses an estimation for the model such that $\|\widehat{\Theta}^{(r)}_j - \Theta_j\| \leq (\frac{1}{2} - \alpha^{(r)})\Delta_{\min}$, for all $j \in [K]$ with $\alpha^{(r)} \in \mathbb{R}$. Moreover, according to Theorem \ref{theorem:convergence}, we have
\begin{align*}
 \| \widehat{\Theta}^{(r+1)}_j -\Theta_j \| &\leq \frac{1}{2}\|  \widehat{\Theta}^{(r)}_j -\Theta_j\| + \Bar{c}_0 \times \frac{1}{\sqrt{\sum_{i \in \widehat{\mathcal{C}}^{(r)}_j} N_i}} \\
    &+\bar{c}_1\Delta_{\max} \sum_{i \in [M] }\sum_{t=0}^{T-1} \exp\left(-\bar{c}_2 N_in_x \left(\frac{\alpha^{(r)} \rho^{(i)} \|\Sigma^{(i)}_t\|}{\rho^{(i)}\|\Sigma^{(i)}_t\| + \sqrt{n_x} }\right)^2\right),
\end{align*}
where by using Assumption \ref{assumption:Ni_nx} and $0<\alpha^{(0)}<\frac{1}{2}$ we can guarantee that $\| \widehat{\Theta}^{(r+1)}_j -\Theta_j \| \leq \|  \widehat{\Theta}^{(r)}_j -\Theta_j\|$ for any $r \in [R]$. This implies that $\alpha^{(r+1)}\geq \alpha^{(r)}$, for any $r \in [R]$.  First, we aim to show that after a small number of iterations, we obtain a sufficiently large value of $\alpha^{(r)} \geq \frac{1}{4}$. To do so, let
\begin{align}\label{eq:epsilon}
    \epsilon_r &:= \Bar{c}_0 \times \frac{1}{\sqrt{\sum_{i \in \widehat{\mathcal{C}}^{(r)}_j} N_i}}  
    +\bar{c}_1\Delta_{\max} \sum_{i \in [M] }\sum_{t=0}^{T-1} \exp\left(-\bar{c}_2 N_in_x \left(\frac{\alpha^{(r)} \rho^{(i)} \|\Sigma^{(i)}_t\|}{\rho^{(i)}\|\Sigma^{(i)}_t\| + \sqrt{n_x} }\right)^2\right),
\end{align}   
be the error at iteration $r$, and note that $\epsilon_{r+1} \leq \epsilon_{r}$ for any $r \in  [R]$ since $\alpha^{(r+1)} \geq \alpha^{(r)}$. Then, after $R^{\prime}$ iterations of Algorithm \ref{algorithm:clustered_fedsysid}, we obtain
\begin{align*}
    \| \widehat{\Theta}^{(R')}_j -\Theta_j \| 
    &\leq (1-\mu_j)^{R'}\left(\frac{1}{2} - \alpha^{(0)}\right)\Delta_{\min} + 2\epsilon_0
\end{align*}
for $R' \geq 2$. Therefore, we need to guarantee that after $R'\geq 2$ parallel iterations, the right hand side of the above expression is upper bounded by  $\frac{1}{4}\Delta_{\min}$. For the first term, since $0<\alpha^{(0)}<\frac{1}{2}$, it suffices to show that $(\frac{1}{2})^{R'} \leq \frac{1}{4}$, which is satisfied for any $R'\geq 2$. On the other hand, $ 2\epsilon_0 \leq \frac{1}{8}\Delta_{\min}$ follows directly from the minimum separation condition of Assumption \ref{assumption:Ni_nx}. Therefore, we have $\| \widehat{\Theta}^{(r)}_j -\Theta_j \| \leq \frac{1}{4}\Delta_{\min}$, for any $r \geq R'$. Then, after $R^{\prime \prime} \geq R^{\prime}$, we have 
 \begin{align*}
    \| \widehat{\Theta}^{(R^{\prime \prime})}_j -\Theta_j \|  \leq \left(\frac{1}{2}\right)^{R^{\prime \prime}}\frac{\Delta_{\min}}{4} + 2\epsilon_0 
 \end{align*}
 which implies $\| \widehat{\Theta}^{(R)}_j -\Theta_j \| \leq \epsilon$ after $R = R^{\prime} + R^{\prime  \prime} \geq 2+\log(\frac{\Delta_{\min}}{4\epsilon})$, with $\epsilon$ as defined in \eqref{eq:epsilon_corollary}.

\end{document}